\def\Alb{{\rm Alb}}
\def\dim{{\rm dim}}
\def\et{\acute{\rm e}{\rm t}}
\def\hra{\hookrightarrow}
\def\Ker{{\rm Ker}}
\def\lra{\longrightarrow}
\def\pr{{\rm pr}}
\def\Ra{\Rightarrow}
\def\ra{\longrightarrow}
\def\rk{{\rm rk}}
\def\skpr{\ltimes}
\def\tor{{\rm tor}}
\documentstyle[12pt,fullpage,pb-diagram,amssymb]{article}

\newtheorem{claim}{Claim}
\newtheorem{lemma}{Lemma}

\newtheorem{cor}{Corollary}
\title{A connectivity lemma for the Albanese map}
\author{M. Rovinsky\thanks{The author was supported 
           by Paul and Gabriella Rosenbaum fellowship} \\
           M.I.T.}
\date{\today}
\begin{document}
\maketitle
\begin{abstract}
I prove that for any complex projective variety $X$ and a 
sufficiently large integer $N$ all the fibers of Albanese map of the 
$N$-th configuration space of $X$ are dominated by smooth connected 
projective varieties with vanishing ${\rm H}^1$. 

This result reduces certain problems on 0-cycles on smooth varieties 
to those on smooth surfaces with vanishing ${\rm H}^1$'s. 
\end{abstract}

\vfill

Let $X$ be a smooth projective complex variety. 
After a choice of a point on $X$ one defines the Albanese map to an 
abelian variety $\Alb X$ (called the {\it Albanese variety} of $X$) 
by the following universality property: any morphism from $X$ to an 
abelian variety $B$, which sends the marked point to the identity 
element of $B$ factors uniquely through a homomorphism $\Alb X\lra B$. 

It is well-known, that in general the fibers of the Albanese map are 
not connected (e.g., for generic hyperplane sections of projective 
embeddings of products of two-dimensional abelian varieties with 
the projective line). However, Lemma \ref{con} below shows that the 
fibers are connected for ``sufficiently iterated'' Albanese map. 

Moreover, one has the following
\begin{claim}
\label{top-red}
Let $X$ be a smooth irreducible projective variety over an 
algebraically closed field of characteristic zero. Then the set of points 
of its Albanese variety is canonically identified with the set of classes 
of degree-zero 0-cycles on $X$ parameterized by smooth projective 
surfaces with vanishing ${\rm H}^1$'s.
\end{claim}

{\it Proof}  is divided into two lemmas, where we keep the following
notations:\\ 
$p: X\ra \Alb X$ is an Albanese map;
$p_N:X^N \ra \Alb X,~~(x_1,...,x_N)\mapsto p(x_1)+...+p(x_N)$.

\begin{lemma}
\label{con}
For $N$ sufficiently large all fibers of $p_N$ are connected.
\end{lemma}

{\it Proof.} Assume that $p_N$ and $p_k$ are surjective, but $p_N$ has
non-connected fibers. 

Then the projection $pr_2$ on $X^k$ of the fiber of $p_{N+k}$ over a point 
$a\in\Alb X$
$$p_{N+k}^{-1}(a)\!=\!
\{(x_1,...,x_{N+k})\!\in \!X^{N+k}\mid p_N(x_1,...,x_N)\!=\!
a-p_k(x_{N+1},...,x_{N+k})\}\stackrel{pr_2}{\lra}X^k$$ 
is also surjective.  

The fiber of $pr_2$ 
over $(x_1,...,x_k)\in X^k$ is $p_N^{-1}(a-p_k(x_1,...,x_k))$.
Since $p_N: X^N\ra \Alb X$ induces  surjection of the fundamental groups,
${\bf Spec}(p_{N\ast}({\cal O}_{X^N}))$ cannot be $\et {\rm ale}$ over 
$\Alb X$  and there is a point of $X^k$ with less quantity of connected 
components of the fiber. Thus any fiber of $p_{N+k}$ has less amount of 
connected components than the generic fiber of $p_N$. \hfill
$\Box$

\begin{lemma}
For any positive integer $N$ there is a complete intersection $Y_N$ of 
dimension 2 with a free action of the symmetric group $\Sigma _N$, 
constructed in \cite{Serre}. 
Let $N$ be as above, 
$\hat{p}_N: (X^N\times Y_N)/\Sigma _N\ra \Alb X,
~~(x_1,...,x_N;y)\mapsto p(x_1)+...+p(x_N)$. 
Then ${\rm H}^1(\hat{p}^{-1}_N(a))=0$ for 
any point $a$ of $\Alb X$.
\end{lemma}

{\it Proof}. One has the Cartesian square
$$
\begin{diagram}
\node{}\node{p_{N+k}^{-1}(a)}\arrow{sw,l}{pr_1}\arrow{se,l}{pr_2}
                         \arrow[2]{e,t,J}{i}
\node{}\node{X^{N+k}}\\
\node{X^{N}}\arrow{se,r}{a-p_N}
\node{\Diamond}\node{X^k}\arrow{sw,r}{p_k}\node{}\node{}\\
\node{}\node{\Alb X}\node{}\node{}\node{}
\end{diagram}
$$
which gives 
$R^q\pr_{1\ast}({\Bbb Q})=(a-p_N)^{\ast}(R^q p_{k\ast}({\Bbb Q}))$. 

Assume that $k$ is large enough and consider the following 
two Leray spectral sequences 
\[
E^{p,q}_2(k,\!N)={\rm H}^p(X^N,R^q\pr_{1\ast}({\Bbb Q}))\Ra
{\rm H}^{p+q}(p_{N+k}^{-1}(a),{\Bbb Q}),
\]
\[
E^{p,q}_2(k)={\rm H}^p(\Alb X,R^qp_{k\ast}({\Bbb Q}))\Ra
{\rm H}^{p+q}(X^k,{\Bbb Q}).
\]
Clearly, $(a-p_N)^{\ast}:E^{p,q}_2(k)\ra E^{p,q}_2(k,N)$ induces an 
isomorphism on $E^{0,q}_2$ and commutes with differentials. Due to the 
Hodge decomposition $p^{\ast}_k:E_2^{p,0}(k)\hra {\rm H}^p(X^k,{\Bbb Q})$. 
So $E_2^{p,0}(k)$ coincides with $E_{\infty}^{p,0}(k)$ and
$d_2^{p-2,1}(k)=0$.

From the commutative diagram
$$
\begin{diagram}
\node{E^{0,1}_2(k)} \arrow[2]{e,t}{d_2(k)=0} \arrow{s,l}{\cong}
\node[2]{E^{2,0}_2(k)} \arrow{s} \\
\node{E^{0,1}_2(k,N)} \arrow[2]{e,t}{d_2(k,N)}
\node[2]{E^{2,0}_2(k,N)}
\end{diagram}
$$
we get $d_2^{0,1}(k,N)=0$, and 
$E_{\infty}^{0,1}(k,N)=E_2^{0,1}(k,N)\cong E_2^{0,1}(k)\cong 
{\rm H}^1(X^k,{\Bbb Q})/E_2^{1,0}(k)$. 

Thus, $\dim {\rm H}^1(p_{N+k}^{-1}(a),{\Bbb Q})\!=\!\dim E^{1,0}_2(k,N)+
\dim ({\rm H}^1(X^k,{\Bbb Q})/E^{1,0}_2(k))$. 

Finally, $\dim {\rm H}^1(p_N^{-1}(a),{\Bbb Q})=(N-1)\rk {\rm H}^1(X)$.

On the other hand, there is a commutative diagram
\[
\begin{diagram}
\node{{\rm H}^1(X^{N+k},{\Bbb Q})} \arrow[2]{e,t}{i^{\ast}}
\node[2]{{\rm H}^1(p_{N+k}^{-1}(a),{\Bbb Q})}\\
\node{{\rm H}^1(X^N,{\Bbb Q})} \arrow{n} \arrow{ene,r}{\pr _1^{\ast}}
\end{diagram}
\]
where $\pr _1^{\ast}$ is injective as $E_2^{1,0}$ coincides with 
$E_{\infty}^{1,0}$ for any spectral sequence in the first quadrant, and 
$\dim {\rm H}^1(X^k,{\Bbb Q})\geq \dim \Ker i^{\ast}$. 
But $\Ker i^{\ast}$ is a $\Sigma _{N+k}$-invariant subspace 
and must have dimension at least $\rk {\rm H}^1(X)$. So it is diagonal 
of ${\rm H}^1(X^{N+k})={\rm H}^1(X)^{N+k}$. 
From the exact sequence
\[
0\ra {\rm H}^1(X)\stackrel{\Delta}{\longrightarrow}
{\rm H}^1(X^N)\ra {\rm H}^1(p_N^{-1}(a))\ra 0
\] 
one gets ${\rm H}^1(p_N^{-1}(a))^{\Sigma _N}=0$, or 
${\rm H}^1(\hat{p}_N^{-1}(a))=0$. \hfill
$\Box$

Note that almost all fibers of $\hat{p}_N$ are smooth. For any $N$, any
point of the critical locus $\Theta _N\subset\Alb X$ of the map $p_N$ can
be translated out of $\Theta _{N+k}$ for some $k$. Then the Lefschetz
hyperplane section theorem gives the rest of the proof of the claim. 
\hfill 
$\Box$

\vspace{7mm}

{\sc Remark.} For locally trivial bundles over Albanese variety one could 
proceed as follows. 
If a subgroup $\Gamma$ of $\pi _1(X)$ is defined by the canonical exact
sequence $$1\ra\Gamma\ra\pi _1(X)\ra{\rm H}_1(X)/\tor\ra 0$$ 
and $\Gamma _N$ is also defined by again the canonical exact
sequence 
$$1\ra\Gamma _N\ra\pi _1(X)^N\skpr\Sigma _N\ra{\rm H}_1(X)/\tor\ra 0,$$
then $\pi _1(\hat{p}_N^{-1}(a))=\Gamma _N$ is the fundamental group of the
fiber of $\hat{p}_N$ and 
$$\Gamma _N=\{(g_1,\dots ,g_N;\sigma)\in\pi _1(X)^N\skpr\Sigma _N
\mid g_1\cdots g_N\in\Gamma\},$$
where $\Sigma _N$ is the symmetric group on $N$ elements.
For $N\geq 3$ 
$$[\Gamma _N,\Gamma _N]=
\{(g_1,\dots ,g_N;\sigma)\in\pi _1(X)^N\skpr A_N\mid 
g_1\cdots g_N\in [\pi _1(X),\pi _1(X)]\},$$
where $A_N$ is the alternating group.
Thus
$${\rm H}_1(\hat{p}_N^{-1}(a))=\Gamma _N/[\Gamma _N,\Gamma _N]\cong 
{\rm H}_1(X)_{\tor}\times {\bf Z}/2{\bf Z}.$$

\end{document}